\documentclass[10pt, a4paper,reqno]{amsart}

\usepackage{amsmath}
\usepackage{amsfonts}
\usepackage{amssymb}
\usepackage{amsthm}
\usepackage{latexsym}
\usepackage[active]{srcltx}
\usepackage{color}
\usepackage[colorlinks,citecolor=blue,linkcolor=blue]{hyperref}

\numberwithin{equation}{section}

\newtheorem{theorem}[equation]{Theorem}
\newtheorem{lemma}[equation]{Lemma}

\theoremstyle{definition}
\newtheorem{definition}[equation]{Definition}

\theoremstyle{remark}
\newtheorem{remark}[equation]{Remark}

\def\kint_#1{\mathchoice%
          {\mathop{\kern 0.2em\vrule width 0.6em height 0.69678ex depth -0.58065ex
                  \kern -0.8em \intop}\nolimits_{\kern -0.4em#1}}%
          {\mathop{\kern 0.1em\vrule width 0.5em height 0.69678ex depth -0.60387ex
                  \kern -0.6em \intop}\nolimits_{#1}}%
          {\mathop{\kern 0.1em\vrule width 0.5em height 0.69678ex depth -0.60387ex
                  \kern -0.6em \intop}\nolimits_{#1}}%
          {\mathop{\kern 0.1em\vrule width 0.5em height 0.69678ex depth -0.60387ex
                  \kern -0.6em \intop}\nolimits_{#1}}}
\def\vintslides_#1{\mathchoice%
          {\mathop{\kern 0.1em\vrule width 0.5em height 0.697ex depth -0.581ex
                  \kern -0.6em \intop}\nolimits_{\kern -0.4em#1}}%
          {\mathop{\kern 0.1em\vrule width 0.3em height 0.697ex depth -0.604ex
                  \kern -0.4em \intop}\nolimits_{#1}}%
          {\mathop{\kern 0.1em\vrule width 0.3em height 0.697ex depth -0.604ex
                  \kern -0.4em \intop}\nolimits_{#1}}%
          {\mathop{\kern 0.1em\vrule width 0.3em height 0.697ex depth -0.604ex
                  \kern -0.4em \intop}\nolimits_{#1}}}

\newcommand{\R}{\mathbb{R}}
\newcommand{\Rn}{\mathbb{R}^d}
\newcommand{\Z}{\mathbb{Z}}

\renewcommand{\l}{\left}
\renewcommand{\r}{\right}

\def\XXint#1#2#3{{\setbox0=\hbox{$#1{#2#3}{\int}$}
\vcenter{\hbox{$#2#3$}}\kern-.5\wd0}}

\newenvironment{list2}{
  \begin{list}{$\bullet$}{%
      \setlength{\itemsep}{0in}
      \setlength{\parsep}{0in} \setlength{\parskip}{0in}
      \setlength{\topsep}{0in} \setlength{\partopsep}{0in} 
      \setlength{\leftmargin}{0.2in}}}{\end{list}}

\title[Positivity of the fundamental solution]{Positivity of the fundamental solution for fractional diffusion and wave equations}
\author{Jukka Kemppainen}

\begin{document}



\begin{abstract}
We study the question of positivity of the fundamental solution for fractional diffusion and wave equations of the form, which may be of fractional order both in space and time. We give a complete characterization for the positivity of the fundamental solution in terms of the order of the time derivative \(\alpha\in(0,2)\), the order of the spatial derivative \(\beta\in (0,2]\) and the spatial dimension \(d\). It turns out that the fundamental solution fails to be positive for all \(\alpha\in (1,2)\) and either \(\beta\in (0,2]\) and \(d\ge 2\) or \(\beta<\alpha\) and \(d=1\), whereas in the other cases it remains positive. 
The proof is based on delicate properties of the Fox H-functions and the Mittag-Leffler functions.
\end{abstract}

\maketitle

\section{Introduction}

We study the positivity of the fundamental solution for the diffusion-wave equation
\begin{equation}\label{prob}
\partial_t^{\alpha} u(t,x)+(-\Delta)^{\beta/2} u(t,x)=0\quad 
\mathrm{in}\quad\R_+ \times \Rn,\ \  0 < \alpha< 2,\ \ 0<\beta\le 2,
\end{equation} 
where $\partial_t^\alpha$ denotes the Caputo fractional derivative if $1\neq\alpha\in (0,2)$ and $\Delta$ is the Laplace operator. When \(\alpha=1\), \(\partial_t^\alpha\) denotes the usual time derivative \(\partial_t\). Since the fractional derivatives are nonlocal, the equation may be nonlocal both in space and time in which case we call it {\em a fully nonlocal diffusion-wave equation}.

We present a unified approach in terms of the Fox $H$-functions and the Mittag-Leffler functions to study the positivity on the whole range \(\alpha\in (0,2),\beta\in(0,2]\) in all spatial dimensions \(d\ge 1\). Large part of the results we present here are already known in the literature but a unified presentation is missing. In particular, large part of the results in the literature concentrate on the {\em subdiffusive region} \(0<\alpha<1\) and to the case \(\beta=2\). The case \(\alpha\le 1\) is well-studied, since it is known that the fundamental solution induces a probability density in \(\Rn\), see e.g.~\cite{KempSiljZach17}. Alternatively, one may use the nice subordination principle to cover the region \(0<\alpha<1\) and \(0<\beta<2\), see e.g.~\cite{Bazh18,Bazh19,Prus93,Tayl10}. For \(\alpha>1\) there are fewer works. In the literature it is already observed that the fundamental solution may not be positive in the multidimensional case, if \(\alpha>1\). But to the best of author's knowledge the results covering the whole range \((\alpha,\beta,d)\in (1,2)\times(0,2]\times\mathbb{Z}_+\) seem to be missing. The celebrated paper of Schneider and Wyss~\cite{SchnWyss89} covers the whole range \(0<\alpha<2\) in the case \(\beta=2\), and it is shown that the fundamental solution changes sign in dimensions \(d\ge 2\), if \(\alpha>1\). See also~\cite{Hany02ii}. The one-dimensional case in the region \(\alpha>1\) and \(\beta=2\) is studied in~\cite{MainParaGore07}, where it is proved that the fundamental solution is positive. For the fully nonlocal diffusion-wave equation the failure of positivity is proved in ~\cite{Hany02} for \(1<\alpha\le\beta<2\) and \(d=3\). The case \(\alpha=\beta\in(1,2)\) in dimensions \(d=1,2,3\) is studied in~\cite{Luch14}, where it is proved that the fundamental solution is a probability density function in dimension \(d=1\) and changes sign in dimension \(d=3\).

In our analysis we will utilize the fine properties of the Fox $H$-functions including simplification rules, the series representations, the asymptotic behavior and the integral transform formulas. The Fox $H$-functions generate a rich family of elementary and special functions, which is a proper framework for studying various properties of the model problem~\eqref{prob}. For example, the asymptotic behavior of the fundamental solution is studied in~\cite{KimLim15,KempSiljVergZach14, KempSiljZach17}, which turn out to be useful.

The author hopes that knowing the borderlines of positivity of the fundamental solution of the simple model problem~\eqref{prob} gives insight for the probabilistic interpretation of more general problems. There is a strong interplay between the stochastic processes and (integro-)PDEs. In particular, the probability density function of the stochastic process at time instant \(t>0\) is a solution of a (integro-)PDE such as~\eqref{prob}.  For this connection we refer to~\cite{MeerSiko12} and references therein. 

The positivity of the fundamental solution or lack of it has also a crucial effect on the properties of the solutions \(u(t,x)\) of~\eqref{prob}. In particular, if the fundamental solution \(G(t,x)\) fails to be positive, the solution operator
\[
S_{\alpha,\beta}^t:\,u_0\mapsto u(x,t)=(G\star u_0)(t,x),\quad u_0(x)=u(x,0),
\] 
where \(\star\) denotes the convolution in the spatial variable \(x\), does not preserve positivity. This implies that in general we do not have nontrivial lower estimates for the norms of the solutions. From physical point of view the lost of positivity means that the equations of the form~\eqref{prob} cannot model diffusion in spatial dimensions greater than one, if \(\alpha>1\). 

\section{Preliminaries and the positivity result} \label{section:preliminaries}


Let us first fix some notations. We denote the space of \(k\)-times 
continuously differentiable functions by \(C^k\) and \(C^0:=C\). 

The Riemann-Liouville fractional integral of order \(\alpha \ge 0\) is defined for $\alpha = 0$ as \(J^0:=I\), where \(I\) denotes the identity operator, and for $\alpha >0$ as
\begin{equation}\label{Riemann-Liouville_int}
J^\alpha 
f(t)=\frac{1}{\Gamma(\alpha)}\int_0^t(t-\tau)^{\alpha-1}f(\tau)\mathrm{d}\tau = (g_\alpha * f)(t),
\end{equation}
where 
\[
g_\alpha(t)=\frac{t^{\alpha-1}}{\Gamma(\alpha)}
\]
is the Riemann-Liouville kernel and $*$ denotes the convolution in time. 

The Caputo fractional derivative of order \(0<\alpha<2\) is defined by
\begin{equation}\label{Caputo_def}
\partial_t^\alpha f(t)=J^{\lceil\alpha\rceil-\alpha}\left(\frac{\mathrm{d}^{\lceil\alpha\rceil}}{\mathrm{d} t^{\lceil\alpha\rceil}}f\right)(t),
\end{equation}
where \(\rceil\alpha\rceil=\min\{k\in\mathbb{Z}\,:\, k\ge \alpha\}\) denotes the least integer greater than or equal to \(\alpha\). Using~\eqref{Riemann-Liouville_int} we can write the definition~\eqref{Caputo_def} in the integral form
\begin{equation}\label{Caputo_def_intform}
\partial_t^\alpha f(t)=\begin{cases}
\frac{1}{\Gamma(1-\alpha)}\int_0^t (t-\tau)^{-\alpha} f'(\tau)\mathrm{d}\tau, &\text{when } 0<\alpha<1,\\
\frac{1}{\Gamma(2-\alpha)}\int_0^t (t-\tau)^{1-\alpha} f''(\tau)\mathrm{d}\tau, &\text{when } 1<\alpha<2,
\end{cases}
\end{equation} 
and, adopting the convention \(J^0:=I\), the definition~\eqref{Caputo_def} reduces to the usual time derivative \(\partial_t\), when \(\alpha=1\).

Let 
\[
\widehat{u}(\xi)=\mathcal{F}(u)(\xi)=(2\pi)^{-d/2}\int_{\Rn}\mathrm{e}^{
-\mathrm{i}x\cdot\xi}f(x)\mathrm{d}x
\] 
and
\[
\mathcal{F}^{-1}(u)(\xi):=\mathcal{F}(u)(-\xi)
\]
denote the Fourier and inverse Fourier transforms of $u$, respectively. We define the fractional Laplacian 
as 
\begin{equation}\label{fract_Laplacian_def}
(-\Delta)^{\beta/2}u(x)=\mathcal{F}^{-1}_{\xi\to x}(|\xi|^\beta
\widehat{u}(\xi)).
\end{equation}
Let \(\mathcal{S}(\Rn)\) denote the Schwarz space, when its dual \(\mathcal{S}'(\Rn)\) is the space of tempered distibutions. Since the Fourier transform is a mapping from \(\mathcal{S}'(\Rn)\) into itself, the fractional Laplacian given by~\eqref{fract_Laplacian_def} defines a mapping from \(\mathcal{S}'(\Rn)\) into itself.

In literature there are different definitions for the fundamental solutions of PDEs, but there seems to be no formal definition of the fundamental solution for the fractional partial differential equations. We will give here a rather weak formulation in terms of being a solution of the equation~\eqref{prob}. The reason is, as it is known and as we shall later see, that the fundamental solution in general is singular at \(x=0\) not only, when \(t=0\), but also for later times, see e.g.~\cite{KempSiljVergZach14,KempSiljZach17,Koch90}. Therefore the fundamental solution is not actually a proper solution of~\eqref{prob}. Since we do not want to dwell into detailed regularity of the solution, we define the concept of the fundamental solution as follows.

\begin{definition}\label{fund_sol_definition}
The function \(\Phi:\R_+\times\Rn\to \R\) is called a {\em fundamental solution} of~\eqref{prob}, if \(\Phi(t,\cdot)\) solves~\eqref{prob} in the sense of \(\mathcal{S}'(\Rn)\) for all \(t>0\) and
\[
\lim_{t\searrow 0} \Phi(t,x)=\delta(x)\quad\text{in  }\mathcal{S}'(\Rn)
\]
together with \(\lim_{t\searrow 0}\partial_t\Phi(t,x)=0\) in \(\mathcal{S}'(\Rn)\), if \(1<\alpha<2\).
\end{definition} 

With this definition at hand, \(\Phi(t,\cdot)\) solves~\eqref{prob} in \(\mathcal{S}'(\Rn)\) if and only if \(\widehat{\Phi}(t,\xi)\) solves
\begin{equation}\label{prob_Fourier}
\partial_t^\alpha \widehat{u}(t,\cdot)+|\cdot|^\beta\widehat{u}(t,\cdot)=0\quad\text{in  } \mathcal{S}'(\Rn)
\end{equation}
for all \(t>0\)~\cite{Tayl96}. It is known that the Mittag-Leffler function \(E_\alpha(-|\cdot|^\beta t^\alpha)\) defined by~\eqref{mittag_asymp1} is a solution of~\eqref{prob_Fourier}~\cite[Chapter 4]{KilbSrivTruj06}. Hence it remains to find a function satisfying the initial condition(s) and whose Fourier transform is the Mittag-Leffler function. The formula for the fundamental solution of~\eqref{prob} given in terms of the Fox $H$-function can be found from the literature, see e.g.~\cite{Duan05,KempSiljZach17,KimLim15}. We will show that the fundamental solution can be written in a form
\begin{equation}\label{fund_solution}
G_{\alpha,\beta,d}(t,x)=\pi^{-d/2}|x|^{-d}H_{23}^{21}\left( 2^{-\beta} 
t^{-\alpha}|x|^{\beta} \big| \begin{smallmatrix} (1,1),& 
(1,\alpha)& \\ (d/2,\beta/2), 
&(1,1), 
&(1,\beta/2) 
\end{smallmatrix} \right),
\end{equation}
where \(H_{23}^{21}\) is a Fox $H$-function whose precise definition is given in Appendix~\ref{Appendix}. We will drop the subindices from $G$ in~\eqref{fund_solution} and simply denote \(G:=G_{\alpha,\beta,d}\) if there is no danger of confusion. We will prove  

\begin{theorem}\label{positivity_theorem}
The fundamental solution of the problem~\eqref{prob} is given by~\eqref{fund_solution}. It
\begin{itemize}
\item[(a)] is positive, if either \(\alpha\in (0,1],\beta\in (0,2]\) and \(d\ge 1\), or \(\alpha\in (1,2)\), \(\beta\in[\alpha,2]\) and \(d=1\);
\item[(b)] changes sign in the following cases of the parameters:
\begin{itemize}
\item[(i)] \(d\ge 2\), \(\alpha\in (1,2)\) and \(\beta\in(0,2]\);
\item[(ii)] \(d=1\), \(\alpha\in (1,2)\) and \(\beta<\alpha\).
\end{itemize}
\end{itemize}
\end{theorem}

Note that, when \(\beta=2\), we can use the property \((iii)\) of Lemma~\ref{Fox_properties} to simplify the expression of the fundamental solution into the form
\begin{equation}\label{fundsol_beta2}
G_{\alpha,2,d}(t,x)=\pi^{-d/2}|x|^{-d}H_{12}^{20}\left( \frac14 
t^{-\alpha}|x|^{2} \big| \begin{smallmatrix} 
(1,\alpha)& \\ (d/2,1), 
&(1,1) 
\end{smallmatrix} \right),
\end{equation}
which is the form given in~\cite{KempSiljVergZach14}, see also~\cite{EideKoch04}. If we further specialize \(\alpha=1\) in~\eqref{fundsol_beta2} and use the symmetry property of the $H$-function \(H^{20}_{12}\) for the parameters appearing in \(A(s)\) of~\eqref{mellin_transform}  together with the property $(iii)$ of Lemma~\ref{Fox_properties}, we obtain the representation
\[
G_{1,2,d}(t,x)=\pi^{-d/2}|x|^{-d}H^{10}_{01}\left( \frac{|x|^2}{4t}\big| \begin{smallmatrix} -\\ (d/2,1)\end{smallmatrix} \right).
\]
Further, using the property $(v)$ with \(\sigma=d/2\) and the formula~\eqref{exp_fox}, the representation~\eqref{fund_solution} reduces to the heat kernel
\[
G_{1,2,d}(t,x)=\frac{1}{(4\pi t)^{d/2}}\exp\left(-\frac{|x|^2}{4t}\right).
\]
\section{Proof of Theorem~\ref{positivity_theorem}}

\subsection{Asymptotic behavior and the Fourier transform of the fundamental solution}

For the moment we assume that \(G\) is given by~\eqref{fund_solution} and call it fundamental solution although we have not proved that \(G\) satisfies the condition(s) of Definition~\ref{fund_sol_definition}. We start with the asymptotic behavior of \(G\). We have the following result

\begin{lemma}\label{fund_sol_asympt}

Let $d \in \Z_+$, $0<\alpha<2$ and $0< \beta \le 2$. Denote \(R:=|x|^\beta 
t^{-\alpha}\). Then the function $G$ given by the formula~\eqref{fund_solution} has the following asymptotic behavior:
\begin{itemize}

\item[(i)] If \(R\leq 1\), then
\[
|G(t,x)|\lesssim \begin{cases}
             t^{-\alpha d/\beta}, &\textrm{if $\alpha=1$, or $\beta > d$ and 
$0<\alpha<2$},\\
t^{-\alpha}(|\log(|x|^\beta t^{-\alpha})|+1), &\textrm{if $\beta=d$ and  
$0<\alpha <2$},\alpha\neq 1,\\
t^{-\alpha}|x|^{-d+\beta} &\textrm{if $0<\beta < d$ and $0<\alpha<2$}.
            \end{cases}
\]

\item[(ii)] If \(R\ge 1\), then 
\[
|G(t,x)|\lesssim t^{\alpha}|x|^{-d-\beta},\quad \textrm{if \(\beta<2\)}.
\]
In the special case $\beta=2$ there holds
\[
|G(t,x)|\lesssim t^{-\alpha d/2}\exp\left(-cR^{\frac{1}{2-\alpha}}\right)
\]
for some positive constant \(c\).

\end{itemize}

\end{lemma}

\begin{proof}
The proof is given in~\cite{KempSiljZach17} in the case \(\alpha\le 1\). A closer look at the proof reveals that the role of \(\alpha\) is merely a parameter as long as \(\alpha<2\). The asymptotics follow from the results and discussion of Appendix. For the convenience of the reader we sketch here the idea of the proof.

The asymptotics in \((i)\) follow from Cauchy's residue theorem and the 
Mellin transform
\begin{align*}
\mathcal{H}^{21}_{23}(s):=& \mathcal{M}\left(H^{21}_{23}\big[ z \left[ 
\begin{smallmatrix}
       (1,1),&(1,\alpha) &\\
(\frac{d}{2},\frac{\beta}{2}), &(1,1), &(1,\frac{\beta}{2})
      \end{smallmatrix}
\right]\right)(s)\\
=&\frac{\Gamma(\frac{d}{2}+\frac{\beta}{2}s)\Gamma(1+s)\Gamma(-s)}{
\Gamma(1+\alpha s)\Gamma(-\frac{\beta}{2}s)}
\end{align*}
of the $H$-function appearing in~\eqref{fund_solution}, see also formulae~\eqref{mellin_barnes_integral} and~\eqref{mellin_transform}. Since the Gamma function has simple poles at non-positive integers, the leading term of the series expansion is determined by the poles of the function \(s\mapsto \mathcal{H}^{21}_{23}(s)z^{-s}\). The different cases of (i) are consequences of the different possibilities for the leading term. 

For example, if \(\alpha=1\) or \(\beta>d\), the leading term is given by the pole of the Gamma function \(\Gamma(\frac{d}{2}+\frac{\beta}{2}s)\) at \(s=-\frac{d}{\beta}\), since in the preceeding case \(\Gamma(1+\alpha s)\) in the denominator cancels \(\Gamma(1+s)\) in the numerator, whereas in the latter case the pole of the Gamma function \(\Gamma(\frac{d}{2}+\frac{\beta}{2}s)\) at \(s=-\frac{d}{\beta}\) gives the dominating power \(z^{d/\beta}\). Then
\[
G(t,x)\sim\pi^{-d/2}|x|^{-d}\mathrm{Res}_{s=-\frac{d}{\beta}}\mathcal{H}_{23}^{21}(s)(2^{-\beta} R)^{-s}=\pi^{-d/2}\frac{2^{1-d}\Gamma(1-\frac{d}{\beta})\Gamma(\frac{d}{\beta})}{\beta\Gamma(1-\frac{\alpha d}{\beta})\Gamma(\frac{d}{2})}t^{-\alpha d/\beta},
\] 
as \(R=|x|^\beta t^{-\alpha}\to 0\). The other cases can be treated in a similar fashion.

The asymptotics in \((ii)\) follow from Theorem~\ref{Braa36} for \(\beta\neq 2\). When \(\beta=2\), then by Theorem~\ref{thm_exp_asymptotics} the $H$-function \(H^{20}_{12}\) in~\eqref{fundsol_beta2} decays exponentially at infinity \(R\to\infty\). For details we refer to~\cite{Braa36} and~\cite[Section 1.7]{KilbSaig04}, see also~\cite{EideKoch04,KempSiljVergZach14}.
\end{proof}

\begin{remark}
\begin{itemize}
\item[(i)] Note that in~\cite{KempSiljZach17} we had sharp asymptotics due to the positivity of the fundamental solution. When \(\alpha>1\), the fundamental solution may no longer be positive, whence in general we have only upper estimates in Lemma~\ref{fund_sol_asympt}. But if we replace the conditions \(R\le 1\) and \(R\ge 1\) with the limiting forms \(R\to 0\) and \(R\to\infty\), we obtain sharp asymptotical behavior. We will use this observation later, when we prove that the fundamental solution does not preserve positivity.

\item[(ii)] Since the Fox H-function \(H^{21}_{23}(z)\) appearing in~\eqref{fund_solution}  is an analytic function of \(z\neq 0\), \(G(t,x)\) is a smooth function as long as \(x\neq 0\). But in the spatially multidimensional case \(G(t,\cdot)\) has singularity at the origin for all times \(t>0\). Hence the fundamental solution cannot actually be a proper solution of~\eqref{prob}. 
\end{itemize}
\end{remark}

The estimates \((i)\) of Lemma~\ref{fund_sol_asympt} imply that \(G(t,\cdot)\) is locally integrable for any fixed \(t>0\). The estimates \((ii)\) of Lemma~\ref{fund_sol_asympt} in turn imply that \(G(t,\cdot)\) is integrable at infinity. Moreover, since \(G(t,x)\) of~\eqref{fund_solution} is a smooth function on \(t>0\) and \(x\neq 0\), the function \(G(t,\cdot)\) is integrable in \(\Rn\), so, in particular, the Fourier transform of \(G(t,\cdot)\) exists. Let us calculate the Fourier transform.

We use the Fourier transform formula for a radial function \(f(r):=f(x)\) with \(r=|x|\) defined on \(\Rn\)~\cite{Graf04}
\begin{equation}\label{Fourier_trans_radial}
\widehat{f}(\xi)=|\xi|^{1-\frac{d}{2}}\int_0^\infty  r^{\frac{d}{2}}J_{\frac{d}{2}-1}(r|\xi|)f(r)\mathrm{d}r
\end{equation}
and use the property \((vii)\) of Lemma~\ref{Fox_properties} to obtain
\begin{equation}\label{Fourier_trans_fundsol1}
\widehat{G}(t,\xi)=(2\pi)^{-d/2} H^{22}_{34}\left(t^{-\alpha}|\xi|^{-\beta} \big| \begin{smallmatrix} 
(1,\beta/2), & (1,1) & (1,\alpha)& (d/2,\beta/2)\\ (d/2,\beta/2), &(1,1), &(1,\beta/2) & 
\end{smallmatrix}\right).
\end{equation}
Using the properties \((ii)\), \((iii)\) and \((iv)\) of Lemma~\ref{Fox_properties} we can simplify the Fox H-function appearing in~\eqref{Fourier_trans_fundsol1} as follows
\begin{align*}
H^{22}_{34}\left(t^{-\alpha}|\xi|^{-\beta}\right)&:=H^{22}_{34}\left(t^{-\alpha}|\xi|^{-\beta} \big| \begin{smallmatrix} 
(1,\beta/2), & (1,1) & (1,\alpha)& (d/2,\beta/2)\\ (d/2,\beta/2), &(1,1), &(1,\beta/2) & 
\end{smallmatrix}\right)\\
&=H^{21}_{32}\left(t^{-\alpha}|\xi|^{-\beta} \big| \begin{smallmatrix} 
(1,1) & (1,\alpha)& (d/2,\beta/2)\\ (d/2,\beta/2), &(1,1) &
\end{smallmatrix}\right)\\
&=H^{11}_{21}\left(t^{-\alpha}|\xi|^{-\beta} \big| \begin{smallmatrix} 
(1,1), &(1,\alpha)\\ (1,1) &
\end{smallmatrix}\right)\\
&=H^{11}_{12}\left(t^{\alpha}|\xi|^{\beta} \big| \begin{smallmatrix} 
(0,1)& \\ (0,1), &(0,\alpha)
\end{smallmatrix}\right),
\end{align*}
which together with the formula~\eqref{ML_fox} implies that the Fourier transform of function \(G(t,\cdot)\) defined by~\eqref{fund_solution} is a constant multiple of the Mittag-Leffler function,
\begin{equation}\label{fundsol_fourier}
\widehat{G}(t,\xi)=(2\pi)^{-d/2}E_\alpha(-|\xi|^\beta t^{-\alpha}).
\end{equation}

\subsection{The conditions of Definition~\ref{fund_sol_definition}}

Since the Mittag-Leffler function satisfies~\eqref{prob_Fourier}, it is enough to check that the initial conditions in Definition~\ref{fund_sol_definition} are satisfied.

We see from~\eqref{fund_solution} that \(G\) has the following scaling structure
\begin{equation}\label{fund_sol_scaling}
G(t,x)=t^{-\alpha d/\beta} G(1,t^{-\alpha/\beta} x),\quad t>0,\quad x\neq 0,
\end{equation}
which together with \(G(t,\cdot)\in L^1(\Rn)\) implies
\[
\langle G(t,\cdot),\varphi\rangle =\int_{\Rn} G(1,x)\varphi(t^{\alpha/\beta}x)\mathrm{d}x\to \varphi(0),\quad t\to 0,
\]
for all \(\varphi\in\mathcal{S}(\Rn)\). Here we used the dominated convergence theorem and the fact (see~\eqref{mittag_asymp1})
\[
\int_{\Rn} G(1,x)\mathrm{d}x=(2\pi)^{d/2}\widehat{G}(1,0)=E_\alpha(0)=1,
\]
which holds for all \(0<\alpha<2\). Hence the initial condition for \(G\) is satisfied. 

If \(\alpha>1\), we need also the intial condition for the first time derivative, which due to the smoothness of \(G\) is equivalent to
\[
\langle \partial_t G(t,\cdot),\varphi\rangle=\partial_t\langle\mathcal{F}^{-1} \mathcal{F} G(t,\cdot),\varphi\rangle=\partial_t\int_{\Rn}\widehat{G}(t,\xi)\widehat{\varphi}(-\xi)\mathrm{d}\xi.
\]
Using the series expansion~\eqref{mittag_asymp1}  for the Mittag-Leffler function we obtain
\[
-\frac{\alpha}{\Gamma(1+\alpha)}t^{\alpha-1}\int_{\Rn} |\xi|^\beta\widehat{\varphi}(-\xi)\mathrm{d}\xi+\mathcal{O}(t^{2\alpha-1})\to 0,\quad  t\to 0,
\]
which proves that the second initial condition is satisfied, when \(\alpha>1\).

\subsection{Proof of (a)}

{\em Case} \(\alpha\le 1\): We start with the easiest case \(\alpha\le 1\). In this case we are able to use Schoenberg's Theorem~\ref{thm_Schoenberg}. Since the Fourier transform of the fundamental solution is the Mittag-Leffler function, it is enough to show that the function
\begin{equation}\label{apu1}
\phi(R)=E_\alpha(-t^\alpha R^{\beta/2}),\quad R=|\xi|,
\end{equation}
is completely monotone. In~\cite{Poll38} it is shown that the function \(R\mapsto E_\alpha(-R)\) is completely monotone for \(0<\alpha\le 1\). Since \(R^{\beta/2}\ge 0\) for \(R\ge0\) and \(\frac{\mathrm{d}}{\mathrm{d} R} R^{\beta/2}=\frac{\beta}{2}R^{\frac{\beta}{2}-1}\) is clearly a completely monotone function for \(\beta\le 2\), it follows from Theorem~\ref{thm_Miller_Samko} that \(\phi\) given by~\eqref{apu1} is completely monotone for all \(t>0\).

{\em Case \(d=1\) and \(\beta=\alpha>1\)}: This case is considered in~\cite{Luch13}. For reader's convenience we repeat briefly the argument. It turns out that we are able to derive an explicit formula for the fundamental solution in terms of elementary functions. The formula~\eqref{fund_solution} in this case reads
\begin{equation}\label{fund_solution_d1_alphabeta}
G(t,x)=\pi^{-1/2}|x|^{-1}H_{23}^{21}\left( 2^{-\beta} 
t^{-\beta}|x|^{\beta} \big| \begin{smallmatrix} (1,1),& 
(1,\beta)& \\ (d/2,\beta/2), 
&(1,1), 
&(1,\beta/2) 
\end{smallmatrix} \right).
\end{equation}
Using the definition of the $H$-function~\eqref{mellin_barnes_integral} the formula~\eqref{fund_solution_d1_alphabeta} can be written in terms of the Mellin-Barnes integral
\begin{equation}\label{fund_solution_d1_alphabeta2}
G(t,x)=\frac{1}{2\pi\sqrt{\pi}|x|\mathrm{i}}\int_{\mathcal{L}}\frac{\Gamma(\frac12+\frac{\beta}{2}s)\Gamma(1+s)\Gamma(-s)}{\Gamma(1+\beta s)\Gamma(-\frac{\beta}{2}s)}2^{\beta s}\left(\frac{x}{t}\right)^{-\beta s}\mathrm{d} s.
\end{equation}
Since now the parameters \(D\) and \(\delta\) defined as~\eqref{D_def} and~\eqref{delta_def} for the $H$ function $H^{21}_{23}$ are $D=\beta-\alpha=0$ and \(\delta=2^{-\beta}\), it follows from Theorem~\ref{fox_series_expansion} that we can determine the series expansion of~\eqref{fund_solution_d1_alphabeta2} for all \(|x|<t\). But before doing that it is better to simplify the ratio of the Gamma functions.

Indeed, using the Legendre duplication formula
\[
\Gamma(z)\Gamma(z+\frac12)=2^{1-2z}\sqrt{\pi}\Gamma(2z)
\]
with \(z=\frac12+\frac{\beta}{2}s\), the formula~\eqref{fund_solution_d1_alphabeta2} can be written in a form
\begin{equation}\label{fund_solution_d1_alphabeta3}
G(t,x)=\frac{1}{|x|}\frac{1}{2\pi\mathrm{i}}\int_{\mathcal{L}}\frac{\Gamma(1+s)\Gamma(-s)}{\Gamma(1+\frac{\beta}{2} s)\Gamma(-\frac{\beta}{2}s)}\left(\frac{x}{t}\right)^{-\beta s}\mathrm{d} s.
\end{equation}
The series expansion is now determined by~\eqref{Fox_series} of Lemma~\ref{Fox_properties} by calculating the residues of the Gamma function \(\Gamma(1+s)\) at \(s=-k\), \(k=1,2,\dots\). Since 
\[
\mathrm{Res}_{s=-k}\Gamma(1+s)=\frac{(-1)^{k-1}}{\Gamma(k)}, 
\]
it follows from Cauchy's residue theorem that
\[
G(t,x)=\frac{1}{|x|}\sum_{k=1}^\infty\frac{(-1)^{k-1}}{\Gamma(1-\frac{\beta}{2}k)\Gamma(\frac{\beta}{2}k)}\left(\frac{|x|}{t}\right)^{\beta k},\quad |x|<t.
\]
Then the Euler's reflection formula
\[
\Gamma(1-z)\Gamma(z)=\frac{\pi}{\sin(\pi z)}
\]
with \(z=\frac{\beta}{2}s\) implies the series expansion
\begin{equation}\label{fund_solution_d1_alphabeta_series}
G(t,x)=-\frac{1}{\pi|x|}\sum_{k=1}^\infty\sin(\beta\pi k/2)\left(-\frac{|x|}{t}\right)^{\beta k}.
\end{equation}
Finally, the summation formula
\[
\sum_{k=1}^\infty r^k\sin(ka)=\mathrm{Im}\left(r^k\mathrm{e}^{\mathrm{i}ka}\right)=\frac{r\sin a}{1-2r\cos a+r^2},\quad a\in\mathbb{R},\quad |r|<1,
\]
implies that the fundamental solution can be represented in a nice form containing only elementary functions,
\begin{equation}\label{fund_solution_d1_alphabeta_final}
G(x,t)=\frac{1}{\pi}\frac{|x|^{\beta-1}t^\beta\sin(\pi\beta/2)}{t^{2\beta}+2|x|^\beta t^\beta\cos(\pi\beta/2)+|x|^{2\beta}}
\end{equation}
for all $|x|<t$.

Using the property $(iii)$ of Lemma~\ref{Fox_properties} and proceeding similarly as above allows to derive~\eqref{fund_solution_d1_alphabeta_final} also in the range \(t<|x|\) and analytic continuation implies that~\eqref{fund_solution_d1_alphabeta_final} holds also for \(|x|=t\). 

Since \(\beta<2\), it follows from~\eqref{fund_solution_d1_alphabeta_final} that the fundamental solution is positive.

{\em Case $d=1$ and $1<\alpha\le\beta\le 2$} 
It is remarkable that this case can be reduced to the previous one. This was observed in~\cite{MainLuchPagn01}, where also a skewness parameter appears in the space-fractional derivative. Here we consider only the fractional Laplacian~\eqref{fract_Laplacian_def}. For the result we calculate the Mellin transform of~\eqref{fund_solution}. Due to the scaling 
relation~\eqref{fund_sol_scaling} and to the fact that \(G_{\alpha,\beta,1}(t,\cdot)\) is an even function, it is enough to prove that
\begin{equation}\label{KG_relation}
K_{\alpha,\beta}(x):=G_{\alpha,\beta,1}(x),\quad x>0,
\end{equation}
is positive. We see from~\eqref{fund_solution} that \(K_{\alpha,\beta}\) is of the form
\[
K_{\alpha,\beta}(x)=\frac{1}{\sqrt{\pi}}x^{-1}f(2^{-\beta} x^\beta).
\]
Then a simple change of variables in the definition of the Mellin transform~\eqref{mellin_trans_def} implies
\begin{equation}\label{eq_vali1}
\mathcal{M}(K_{\alpha,\beta})(s)=\frac{1}{\sqrt{\pi}}\frac{2^{s-1}}{\beta}\mathcal{H}^{21}_{23}\left(\frac{s-1}{\beta}\right),
\end{equation}
where \(\mathcal{H}^{21}_{23}\) denotes the Mellin transform of the $H$ function $H^{21}_{23}$ appearing in~\eqref{fund_solution}. Using~\eqref{mellin_transform}, the equation~\eqref{eq_vali1} reads
\begin{equation}\label{eq_vali2}
\mathcal{M}(K_{\alpha,\beta})(s)=\frac{1}{\sqrt{\pi}}\frac{2^{s-1}}{\beta}\frac{\Gamma(\frac12 s)\Gamma(1-\frac{1}{\beta}+\frac{1}{\beta}s)\Gamma(\frac{1}{\beta}-\frac{1}{\beta}s)}{\Gamma(1-\frac{\alpha}{\beta}+\frac{\alpha}{\beta}s)\Gamma(\frac12-\frac12s)}
\end{equation}
Plugging \(\alpha=\beta\) and comparing \(K_{\alpha,\alpha}\) with \(K_{\alpha,\beta}\), we see that
\[
(\mathcal{M}K_{\alpha,\beta})(s)=(\mathcal{M}K_{\alpha,\alpha})(s)\frac{\Gamma(s)}{\Gamma(1-\frac{\alpha}{\beta}+\frac{\alpha}{\beta}s)}.
\]
Therefore \(K_{\alpha,\beta}\) can be expressed as the Mellin convolution
\begin{equation}\label{vali3}
K_{\alpha,\beta}(x)=\int_0^\infty K_{\alpha,\alpha}\left(\frac{x}{y}\right)M_{\alpha/\beta}(y)\frac{\mathrm{d} y}{y},
\end{equation}
where
\begin{equation}\label{M_function1}
M_{\alpha/\beta}(x)=\mathcal{M}^{-1}\left(\frac{\Gamma(s)}{\Gamma(1-\frac{\alpha}{\beta}+\frac{\alpha}{\beta}s)}\right)(x).
\end{equation}
Since \(K_{\alpha,\alpha}\) is positive, it is enough to show that \(M_{\alpha/\beta}\) is positive. The function \(M_{\alpha/\beta}\) may be viewed as the Fox $H$-function
\begin{equation}\label{M_function_fox}
M_{\alpha/\beta}(x)=H^{10}_{11}\left(x \big|\begin{smallmatrix} (1-\frac{\alpha}{\beta},\frac{\alpha}{\beta})\\ (0,1)\end{smallmatrix}\right).
\end{equation}
Since for \(M_{\alpha/\beta}\) the parameters \(a^*\) and \(D\) defined by~\eqref{astar_def} and~\eqref{D_def} are \(a^*=D=1-\frac{\alpha}{\beta}>0\), then by Theorem~\ref{fox_series_expansion} we can obtain the series expansion by Cauchy's residue theorem as follows
\begin{equation}\label{M_function_series}
M_{\alpha/\beta}(x)=\sum_{k=0}^\infty\mathrm{Res}_{s=-k}\frac{\Gamma(s)}{\Gamma(\frac{\alpha}{\beta}s+1-\frac{\alpha}{\beta})}x^{-s}=\sum_{k=0}^\infty \frac{(-1)^k}{k!\Gamma(-\frac{\alpha}{\beta} k+1-\frac{\alpha}{\beta})} x^k,
\end{equation}
which can be recognized as the so called {\em Mainardi function} 
or a {\em Wright function} \(W_{-\frac{\alpha}{\beta},1-\frac{\alpha}{\beta}}(-x)\)~\cite{MainParaGore07}. One can show by using the properties of the $H$ function in Lemma~\ref{Fox_properties} and the Laplace transform formula for the $H$ function given in Theorem~\ref{thm_H_function_laplace} that the Mainardi function is connected to an \(\alpha/\beta\) stable distribution via the formula~\cite[Formula (A.41)]{MainParaGore07}
\[
\frac{\alpha}{\beta}\int_0^\infty x^{-\frac{\alpha}{\beta}-1}M_{\alpha/\beta}(x^{-\alpha/\beta})\mathrm{e}^{-px}\mathrm{d}x=\mathrm{e}^{-p^{\alpha/\beta}},
\]
which implies positivity of \(M_{\alpha/\beta}\) and therefore the positivity of \(G_{\alpha,\beta,1}(t,x)\) for \(\alpha<\beta\).

{\bf Proof of (b)}. We follow the argument presented in~\cite{SchnWyss89}. We will refer to Theorem~\ref{thm_Bernstein} to conclude the lack of positivity of the fundamental solution. It turns out that instead of a direct attack by calculating the Laplace transform of the fundamental immediately, it is useful to do the same transformations as in~\cite{SchnWyss89} first and then refer to Theorem~\ref{thm_Bernstein}. These transformations done in~\cite{SchnWyss89} enabled the reduction to the Macdonald function on the Laplace transform side. As we show, somewhat surprisingly the same method works also for the more general $H$ functions.

By rotational invariance, it is equivalent to study positivity of the function
\[
\rho_{\alpha,\beta}(t,r)=\frac{2\pi^{d/2}}{\Gamma(\frac{d}{2})}G(t,x)r^{d-1}=\frac{2}{\Gamma(\frac{d}{2})}r^{-1}H_{23}^{21}\left( 2^{-\beta} 
t^{-\alpha}r^{\beta} \big| \begin{smallmatrix} (1,1),& 
(1,\alpha)& \\ (d/2,\beta/2), 
&(1,1), 
&(1,\beta/2) 
\end{smallmatrix} \right)
\]
on \(\mathbb{R}_+\). Further, by substituting \(r=x^{-\alpha/\beta}\) it is equivalent to study positivity of the function 
\[
f_{\alpha,\beta}(t,x)=\frac{\alpha}{\beta}x^{-\frac{\alpha}{\beta}-1}\rho_{\alpha,\beta}(t, x^{-\frac{\alpha}{\beta}}).
\]
Using the property \((vi)\) of Lemma~\ref{Fox_properties}, we can write \(f_{\alpha,\beta}\) in a form
\[
f_{\alpha,\beta}(t,x)=t\psi_{\alpha,\beta}(tx),
\]
where
\begin{equation}\label{psi_def}
\psi_{\alpha,\beta}(x)=\frac{2}{\beta\Gamma(\frac{d}{2})}x^{-1} H^{21}_{23}\left(2^{-\beta/\alpha}x^{-1}\big|\begin{smallmatrix} (1,\frac{1}{\alpha}),& 
(1,1)& \\ (d/2,\frac{\beta}{2\alpha}), 
&(1,\frac{1}{\alpha}), 
&(1,\frac{\beta}{2\alpha}) 
\end{smallmatrix} \right).
\end{equation}
Therefore it is enough to study the positivity of the function \(\psi_{\alpha,\beta}\). Using the properties \((v)\) and \((iv)\) of Lemma~\ref{Fox_properties}, we can write \(\psi_{\alpha,\beta}\) in a form
\begin{equation}\label{psi_def2}
\psi_{\alpha,\beta}(x)=\frac{2^{\frac{\beta}{\alpha}+1}}{\beta\Gamma(\frac{d}{2})}H^{12}_{32}\left(2^{\frac{\beta}{\alpha}}x \big|\begin{smallmatrix}
(1-\frac{d}{2}-\frac{\beta}{2\alpha},\frac{\beta}{2\alpha}), & (-\frac{1}{\alpha},\frac{1}{\alpha}), &(-\frac{\beta}{2\alpha},\frac{\beta}{2\alpha})\\ (-\frac{1}{\alpha},\frac{1}{\alpha}), &(-1,1) &
\end{smallmatrix}\right).
\end{equation}

We note that 
\[
\min_{1\le j\le m}\left\{\frac{\mathrm{Re}(b_j)}{\beta_j}\right\}=-1,
\]
whence the condition on Theorem~\ref{thm_H_function_laplace} is not satisfied. This condition is related to the local integrability of the function \(f\) at zero, which guarantees that the Laplace integral 
\[
\mathcal{L}(f)(p)=\int_0^\infty f(t)\mathrm{e}^{-pt}\mathrm{d}t
\]
makes sense near zero. However, a closer look at the Mellin transform
\begin{align*}
\mathcal{H}^{12}_{32}(s)&:=\mathcal{H}^{12}_{32}\left(s \big|\begin{smallmatrix}
(1-\frac{d}{2}-\frac{\beta}{2\alpha},\frac{\beta}{2\alpha}), & (-\frac{1}{\alpha},\frac{1}{\alpha}), &(-\frac{\beta}{2\alpha},\frac{\beta}{2\alpha})\\ (-\frac{1}{\alpha},\frac{1}{\alpha}), &(-1,1) &
\end{smallmatrix}\right)\\
&=\frac{\Gamma(-\frac{1}{\alpha}+\frac{1}{\alpha}s)\Gamma(\frac{d}{2}+\frac{\beta}{2\alpha}-\frac{\beta}{2\alpha}s)\Gamma(1+\frac{1}{\alpha}-\frac{1}{\alpha}s)}{\Gamma(-\frac{\beta}{2\alpha}+\frac{\beta}{2\alpha}s)\Gamma(-s)}
\end{align*}
of the $H$-function appearing in~\eqref{psi_def2} reveals that the singularity at \(s=1\) is removable. The dominating term in the series expansion~\eqref{Fox_series} is determined by the pole of the Gamma function \(\Gamma(-\frac{1}{\alpha}+\frac{1}{\alpha}s)\) at \(s=1-\alpha\), which implies 
\[
|\psi_{\alpha,\beta}(x)|\lesssim\mathcal{O}(x^{\alpha-1}),\quad x\to 0.
\]
Moreover, since \(a^*=2-\alpha>0\) for the $H$ function in~\eqref{psi_def2}, the function \(\psi_{\alpha,\beta}\) has an algebraic decay at infinity by Theorem~\ref{Braa36}, whence the Laplace transform exists (as it should be clear from the beginning, since \(G(t,\cdot)\) is integrable, which implies the existence of the Laplace transform). This indicates that one should not use the properties of the $H$ functions found from literature as a {\em ''cookbook''} of formulae, but instead carefully investigate the nature of singularities by looking at the Mellin transform \(\mathcal{H}^{mn}_{pq}\) of the $H$ function \(H^{mn}_{pq}\).

We have from Theorem~\ref{thm_H_function_laplace} together with properties $(v)$ and $(iv)$ of Lemma~\ref{Fox_properties} that
\begin{align*}
\mathcal{L}(\psi_{\alpha,\beta})(p)&=\frac{2}{\beta\Gamma(\frac{d}{2})}\cdot\frac{1}{\tilde{p}}H^{13}_{42}\left(\frac{1}{\tilde{p}}\big| \begin{smallmatrix}
(0,1),&(1-\frac{d}{2}-\frac{\beta}{2\alpha},\frac{\beta}{2\alpha}), & (-\frac{1}{\alpha},\frac{1}{\alpha}), &(-\frac{\beta}{2\alpha},\frac{\beta}{2\alpha})\\ (-\frac{1}{\alpha},\frac{1}{\alpha}), &(-1,1) & &
\end{smallmatrix}\right)\Big|_{\tilde{p}=2^{-\frac{\beta}{\alpha}}p}\\
&=\frac{2}{\beta\Gamma(\frac{d}{2})}H^{31}_{24}\left( 2^{-\frac{\beta}{\alpha}}p\big| \begin{smallmatrix}
(1,\frac{1}{\alpha}), &(1,1) & & \\
(0,1), &(\frac{d}{2},\frac{\beta}{2\alpha}), &(1,\frac{1}{\alpha}), & (1,\frac{\beta}{2\alpha})
\end{smallmatrix}\right).
\end{align*}
Further, using the property $(vi)$ of Lemma~\ref{Fox_properties} with \(k=\frac{\beta}{2\alpha}\), we can write the Laplace transform of \(\psi_{\alpha,\beta}\) in a form
\begin{equation}\label{psi_laplace_final}
\widetilde{\psi}_{\alpha,\beta}(p):=\mathcal{L}(\psi_{\alpha,\beta})(p)=\frac{4\alpha}{\beta^2\Gamma(\frac{d}{2})}H^{31}_{24}\left(\frac{p^{\frac{2\alpha}{\beta}}}{4}\Big| \begin{smallmatrix}
(1,\frac{2}{\beta}), &(1,\frac{2\alpha}{\beta}) & & \\
(0,\frac{2\alpha}{\beta}), &(\frac{d}{2},1), &(1,\frac{2}{\beta}), & (1,1)
\end{smallmatrix}\right),
\end{equation}
which for \(\beta=2\), as we shall show after the proof, can be reduced to the modified Bessel function of the second kind, which appears in the proof of positivity by Schneider and Wyss~\cite[Formula B.3]{SchnWyss89}. It remains to study the completely monotonicity of \(\widetilde{\psi}_{\alpha,\beta}\). Now we are in a position to study the different cases of (b) formulated in Theorem~\ref{positivity_theorem}. It turns out that the second derivative of \(\widetilde{\psi}_{\alpha,\beta}\) fails to be non-negative in the aforementioned cases. Using the property $(i)$ of Lemma~\ref{Fox_properties} we have
\[
\widetilde{\psi}^{(2)}_{\alpha,\beta}(p)=\frac{4\alpha}{\beta^2\Gamma(\frac{d}{2})}p^{-2}H^{41}_{35}\left(\frac{p^{\frac{2\alpha}{\beta}}}{4}\Big| \begin{smallmatrix}
(1,\frac{2}{\beta}), &(1,\frac{2\alpha}{\beta}),&(0,\frac{2\alpha}{\beta}) & & \\
(2,\frac{2\alpha}{\beta}), &(0,\frac{2\alpha}{\beta}), &(\frac{d}{2},1), &(1,\frac{2}{\beta}), & (1,1)
\end{smallmatrix}\right)
\]
Further, using symmetry of the parameters in the factors \(A,B,C\) and \(D\) in the definition~\eqref{mellin_barnes_integral} with~\eqref{mellin_transform} and the cancellation property $((ii)$ of Lemma~\ref{Fox_properties}, we have
\begin{equation}\label{psi_laplace_2der}
\widetilde{\psi}^{(2)}_{\alpha,\beta}(p)=\frac{4\alpha}{\beta^2\Gamma(\frac{d}{2})}p^{-2}H^{31}_{24}\left(\frac{p^{\frac{2\alpha}{\beta}}}{4}\Big| \begin{smallmatrix}
(1,\frac{2}{\beta}), &(1,\frac{2\alpha}{\beta}),& & \\
(2,\frac{2\alpha}{\beta}), &(\frac{d}{2},1), &(1,\frac{2}{\beta}), & (1,1)
\end{smallmatrix}\right)
\end{equation}
The behavior of the function $\widetilde{\psi}^{(2)}_{\alpha,\beta}$ depends on the poles of the function
\begin{equation}\label{psi_laplace_2der_mellin}
\begin{split}
f(s)&:=\mathcal{H}^{31}_{24}\left(s\Big| \begin{smallmatrix}
(1,\frac{2}{\beta}), &(1,\frac{2\alpha}{\beta}),& & \\
(2,\frac{2\alpha}{\beta}), &(\frac{d}{2},1), &(1,\frac{2}{\beta}), & (1,1)
\end{smallmatrix}\right)z^{-s}\\
&=\frac{\Gamma(2+\frac{2\alpha}{\beta}s)\Gamma(\frac{d}{2}+s)\Gamma(1+\frac{2}{\beta}s)\Gamma(-\frac{2}{\beta}s)}{\Gamma(1+\frac{2\alpha}{\beta}s)\Gamma(-s)}z^{-s}
\end{split}
\end{equation}
with \(z=\tfrac{1}{4}p^{\frac{2\alpha}{\beta}}\). Since the parameter $D$ defined by~\eqref{D_def} for the $H$ function of~\eqref{psi_laplace_2der} is
\[
D=\sum_{j=1}^q\beta_j-\sum_{i=1}^p\alpha_i=2>0,
\]
the asymptotics near zero is given by formula~\eqref{Fox_series} of Theorem~\ref{fox_series_expansion}. We have different cases.

{\em Case $d\ge 2, \alpha\in (1,2)$ and $\beta\in (0,2]$} Since $d\ge2$, we see from~\eqref{psi_laplace_2der_mellin} that the dominating behavior is determined either by pole at $s=-\frac{d}{2}$ or $s=-\frac{\beta}{2}$, since $\alpha<2$ implies $\frac{\beta}{\alpha}>\frac{\beta}{2}$. We have two subcases. If $\beta<2$ or $d\ge 3$, then the pole $s=-\frac{\beta}{2}$ detemines the dominating term. Calculating the residue we have
\[
\mathrm{Res}_{s=-\frac{\beta}{2}}\frac{\Gamma(2+\frac{2\alpha}{\beta}s)\Gamma(\frac{d}{2}+s)\Gamma(1+\frac{2}{\beta}s)\Gamma(-\frac{2}{\beta}s)}{\Gamma(1+\frac{2\alpha}{\beta}s)\Gamma(-s)}z^{-s}=\frac{\beta\Gamma(2-\alpha)\Gamma(\frac{d-\beta}{2})}{2\Gamma(1-\alpha)\Gamma(\frac{\beta}{2})}z^{\frac{\beta}{2}}.
\] 
Since \(\Gamma(1-\alpha)\) is negative and other terms are positive, \(\widetilde{\psi}^{(2)}_{\alpha,\beta}\) takes negative values near the origin, whence \(\widetilde{\psi}_{\alpha,\beta}\) cannot be completely monotone. Theorem~\ref{thm_Bernstein} then implies that \(G(t,\cdot)\) changes sign.

In the other subcase $d=\beta=2$. In this case the order of the pole at $s=-\frac{\beta}{2}$ is two. Hence the dominating term is given by
\[
\mathrm{Res}_{s=-1}\frac{\Gamma(2+\alpha s)\Gamma(1+s)^2\Gamma(-s)}{\Gamma(1+\alpha s)\Gamma(-s)}z^{-s}=\lim_{s\to -1}\frac{\mathrm{d}}{\mathrm{d}s}\left(\frac{\Gamma(2+\alpha s)\Gamma(2+s)^2}{\Gamma(1+\alpha s)}z^{-s}\right).
\]
Since \(\frac{\mathrm{d}}{\mathrm{d} s} z^{-s}=-z^{-s}\log z\) dominates \(z^{-s}\) near zero, the term
\[
-\frac{\Gamma(2-\alpha)}{\Gamma(1-\alpha)}z\log z
\]
dominates and takes negative values, as \(z\to 0\). Hence \(\widetilde{\psi}_{\alpha,\beta}\) cannot be completely monotone and  \(G(t,\cdot)\) changes sign.

{\em Case $d=1,\alpha\in(1,2)$ and $\beta\le 1$} Again we have two subcases. If $\beta<1$, then \(f\) given by~\eqref{psi_laplace_2der_mellin} has a first order pole at \(s=-\frac{\beta}{2}\) and the dominating term is
\[
\mathrm{Res}_{s=-\frac{\beta}{2}} f(s)=\frac{\beta\Gamma(2-\alpha)\Gamma(\frac{1-\beta}{2})}{2\Gamma(1-\alpha)\Gamma(\frac{\beta}{2})}z^{\frac{\beta}{2}}, \quad z\to 0.
\]
Again we can conclude that \(G(t,\cdot)\) changes sign.

In the other subcase \(d=\beta=1\), whence \(f\) has a second order pole at \(s=-\frac12\). Similar calculation as above for the case \(d=\beta=2\) shows that the term
\[
-\frac{\Gamma(2-\alpha)}{2\Gamma(1-\alpha)\Gamma(\frac12)}z^{\frac12}\log z,\quad z\to 0,
\]
dominates, whence again \(G(t,\cdot)\) changes sign.

{\em Case $d=1,1<\alpha<2$ and $1<\beta<\alpha$} This is the final case. Since $\frac{d}{2}=\frac12<\frac{1}{\alpha}<\frac{\beta}{\alpha}$, the dominating term is given by the first order pole at $s=-\frac12$ of $f$ in~\eqref{psi_laplace_2der_mellin}. Calculating the residue we obtain
\[
\mathrm{Res}_{s=-\frac12}f(s)=\frac{\Gamma(2-\frac{\alpha}{\beta})\Gamma(1-\frac{1}{\beta})\Gamma(\frac{1}{\beta})}{\Gamma(1-\frac{\alpha}{\beta})\Gamma(\frac12)}z^{\frac12}.
\]
Since in this case \(1<\frac{\alpha}{\beta}<2\), we have \(1-\frac{\alpha}{\beta}\in (-1,0)\), whence \(\Gamma(1-\frac{\alpha}{\beta})<0\), while other factors are positive. Hence \(G(t,\cdot)\) changes sign, which finishes the proof of Theorem~\ref{positivity_theorem}.\qed

In the particular case $\beta=2$ the Laplace transform~\eqref{psi_laplace_final} is of the form
\[
\widetilde{\psi}_{\alpha,2}(p)=\mathcal{L}(\psi_{\alpha,2})(p)=\frac{\alpha}{\Gamma(\frac{d}{2})}H^{31}_{24}\left(\frac{p^{\alpha}}{4}\Big| \begin{smallmatrix}
(1,1), &(1,\alpha) & & \\
(0,\alpha), &(\frac{d}{2},1), &(1,1), & (1,1)
\end{smallmatrix}\right),
\]
which can be reduced to
\[
\widetilde{\psi}_{\alpha,2}(p)=\mathcal{L}(\psi_{\alpha,2})(p)=\frac{\alpha}{\Gamma(\frac{d}{2})}H^{30}_{13}\left(\frac{p^{\alpha}}{4}\Big| \begin{smallmatrix}
(1,\alpha) & & \\
(0,\alpha), &(\frac{d}{2},1), &(1,1)
\end{smallmatrix}\right).
\]
On the Mellin transform side we obtain
\[
\begin{split}
\mathcal{H}^{30}_{13}\left(s\Big| \begin{smallmatrix}
(1,\alpha) & & \\
(0,\alpha), &(\frac{d}{2},1), &(1,1)
\end{smallmatrix}\right)&=\frac{\Gamma(\alpha s)\Gamma(\frac{d}{2}+s)\Gamma(1+s)}{\Gamma(1+\alpha s)}=\frac{\Gamma(\frac{d}{2}+s)\Gamma(s)}{\alpha}\\
&=\frac{1}{\alpha}\mathcal{H}^{20}_{02}\left(s\Big| \begin{smallmatrix}
- & \\
(0,1), &(\frac{d}{2},1)
\end{smallmatrix}\right)
\end{split}
\]
by using the definition~\eqref{mellin_transform} and the property \(\Gamma(1+z)=z\Gamma(z)\) of the Gamma function. Hence the Laplace transform of \(\psi_{\alpha,2}\) can be reduced to the modified Macdonald function or the modified Bessel function of the second kind
\begin{equation}\label{psi_laplace_final2}
\widetilde{\psi}_{\alpha,2}(p)=\frac{1}{\Gamma(\frac{d}{2})}H^{20}_{02}\left(\frac{p^\alpha}{4}\Big| \begin{smallmatrix}
- & \\
(0,1), &(\frac{d}{2},1)
\end{smallmatrix}\right)=\frac{2^{1-\frac{d}{2}}}{\Gamma(\frac{d}{2})}p^{\frac{\alpha d}{4}} K_{\frac{d}{2}}(p^{\frac{\alpha}{2}})
\end{equation}
by using~\cite[Formula (2.9.19)]{KilbSaig04}. The formula~\eqref{psi_laplace_final2} coincides with~\cite[Formula B.3]{SchnWyss89}.

\section{Conclusions}

We specified the range of parameters, when the fundamental solution of the fractional diffusion and wave equation is positive and can be regarded as a probability density function. From modelling point of view this means that we were able determine the {\em limits of diffusion}. Indeed, when the fundamental solution lacks positivity, then the solution operator \(S_{\alpha,\beta}^t\) for the Cauchy problem of our model equation~\eqref{prob} does not preserve positivity. Hence, the equation~\eqref{prob} cannot model diffusion, when the fundamental solution changes sign.

As we demonstrated, when \(\alpha\le 1\), there is a nice theory, which allows to conclude positivity by the Fourier transform. The analysis simplifies, since on the Fourier side the more complicated Fox $H$-function appearing in representation of the fundamental solution~\eqref{fund_solution} is replaced by a simpler Mittag-Leffler function. Moreover, the Fourier transform of the fundamental solution \(\widehat{G}(t,\xi)=E_{\alpha}(-|\xi|^{\beta}t^\alpha)\) {\em ''does not see the dimension''} \(d\). Hence, Schoenberg's result Theorem~\ref{thm_Schoenberg} allows to conclude positivity for all dimensions \(d\).

Things change, when \(\alpha\) passes the level \(\alpha=1\). Still, in the multidimensional case \(d\ge 2\) the result follows the endpoint result obtained in the paper of Wyss and Schneider~\cite{SchnWyss89}. The lack of positivity for all \(\beta\in(0,2]\) turned out to be in agreement with that of the fundamental solution
\[
\cos (t\sqrt{-\Delta})\delta(x)
\]
for the wave equation, see~\cite[Section 5.3]{Tayl96}. Somewhat surprisingly, the most interesting phenomenom occured in the one-dimensional case \(d=1\). One could guess that the result would follow the positivity (or non-negativity) of the fundamental solution
\[
\left(\cos t\sqrt{-\partial_x^2}\right)\delta(x)=\frac12\left(\delta(x+t)+\delta(x-t)\right)
\]
of the one-dimensional wave equation. This is the case for all \(\alpha\in (0,2)\) as long as \(\beta\) is kept in its endpoint \(\beta=2\), as it was known already due to the results of~\cite{SchnWyss89}. But when \(\beta\) drops below the level \(\beta=2\), the fundamental solution of our model problem turned out change immediately sign as long as \(\beta<\alpha<2\). The case \(\alpha\le\beta\) was still fine as the positivity was preserved. This case was more delicate, since the basic properties of the $H$ functions did not seem to be enough to prove positivity. We needed a proper integral representation~\eqref{vali3} and the positivity   of \(G_{\beta,\beta,1}(t,\cdot)\), which followed from its simple representation~\eqref{fund_solution_d1_alphabeta_final}.

Although the theory of $H$ functions was useful in obtaining the positivity or lack of it, we remark that there is no direct formula, which could be used to prove positivity or finer properties of the fundamental solutions. In general, simplifications are needed to obtain information on quite complicated $H$ function. In particular, from numerical point of view reduction to other functions are needed, because no packages for the numerical calculations of general $H$ function are available~\cite{Luch17}. From mathematical or physical point of view there is need to further analyze the properties of the solutions for the equations of type~\eqref{prob} or its generalizations. Since~\eqref{prob} {\em ''interpolates''} its elliptic (corresponding to~\eqref{prob} with \(\alpha=0\)), parabolic (corresponding to~\eqref{prob} with \(\alpha=1\)) and hyperbolic (corresponding to~\eqref{prob} with \(\alpha=2\)) counterparts, it is expected that~\eqref{prob} or its generalizations inherit some properties from their elliptic, parabolic and hyperbolic counterparts, see e.g.~\cite{Koch19i,Koch19ii,Luch17}. Here we saw that even positivity is nontrivial.



\section{Appendix}\label{Appendix}

Here we recall some classical results which are needed in the theory.

\subsection{Probability theory}

We need the following results from the probability theory. The concept of a {\em completely monotone function} plays an important role. A real function \(f:\R_+\to\R\) is said to be {\em completely monotone}, if
\begin{equation}\label{cm_cond}
(-1)^nf^{(n)}(t)\ge 0\quad\text{for}\quad 0<t<\infty\quad\text{and}\quad n=0,1,2,\dots.
\end{equation}
Nonnegarive functions whose derivative is completely monotone are called {\em Bernstein functions}. The following result is useful, see~\cite{MillSamk01}.

\begin{theorem}\label{thm_Miller_Samko}
Let \(f\) be a completely monotone function and \(g\) a Bernstein function. Then \(f\circ g\) is completely monotone.
\end{theorem}

Combining the Bochner's Theorem on positive definite functions and Schoenberg's result on connecting the positive definiteness to complete monotonicity for radial functions, there holds~\cite{Scho38}
\begin{theorem}\label{thm_Schoenberg}
A function \(\phi\) is completely monotone on \( (0,\infty)\) if and only \(\Phi(|x|)=\phi(|x|^2)\) is the Fourier transform of a nonnegative Borel measure on \(\Rn\).
\end{theorem}

Another useful theorem is the classical result characterizing the completely monotone functions in terms of the Laplace transform.

\begin{theorem}[Hausdorff-Bernstein-Widder Theorem]\label{thm_Bernstein}
A function \(\phi:[0,\infty)\to\R\) is completely monotone on \( (0,\infty)\) if and only if it is the Laplace transform of a nonnegative Borel measure.
\end{theorem}

\subsection{Fox $H$-functions}

The Fox $H$-functions are special functions of a very general nature and there 
is a natural connection to the fractional calculus, since the
fundamental solutions of the Cauchy problem can be represented in terms of 
them. Since the asymptotic behavior of the Fox \(H\)-functions can be found 
from the literature, the Fox \(H\)-functions have a crucial role also in our 
asymptotic analysis. We collect here some basic facts on these functions.

Let us start with the definition. To simplify the notation we introduce
\[
(a_i,\alpha_i)_{k,p}:=((a_k,\alpha_k),(a_{k+1},\alpha_{k+1}),\dots,(a_p,
\alpha_p))
\]
for the set of parameters appearing in the definition of Fox $H$-functions. In the general case the numbers \(a_i\) and \(\alpha_i\) can be complex but in our considerations they are real. The 
Fox \(H\)-function is defined via a Mellin-Barnes type
integral as
\begin{equation}\label{mellin_barnes_integral}
H^{mn}_{pq} (z):= H^{mn}_{pq}\left( z \big| \begin{smallmatrix}
                           (a_i,\alpha_i)_{1,p}\\
(b_j,\beta_j)_{1,q}                 
                                           \end{smallmatrix}
\right) =\frac{1}{2\pi i}\int_{\mathcal{L}}\mathcal{H}^{mn}_{pq}(s)z^{-s} 
\mathrm{d}s,
\end{equation}
where 
\begin{equation}\label{mellin_transform}
\begin{split}
\mathcal{H}^{mn}_{pq}(s)&:=\mathcal{H}^{mn}_{pq}\left( s \big| \begin{smallmatrix}
                           (a_i,\alpha_i)_{1,p}\\
(b_j,\beta_j)_{1,q}                 
                                           \end{smallmatrix}
\right)\\
&=\frac{\prod_{j=1}^m \Gamma(b_j+\beta_j s)\prod_{i=1}^n 
\Gamma(1-a_i-\alpha_i s)}{\prod_{i=n+1}^p\Gamma(a_i+\alpha_i 
s)\prod_{j=m+1}^q\Gamma(1-b_j-\beta_j s)}\\
&=:\frac{A(s)B(s)}{C(s)D(s)}
\end{split}
\end{equation}
is the Mellin transform 
\begin{equation}\label{mellin_trans_def}
\mathcal{M}(f)(s)=\int_0^\infty t^{s-1} f(t)\mathrm{d}t
\end{equation}
of the Fox $H$-function \(H_{pq}^{mn}\) and 
\(\mathcal{L}\) is the infinite contour in the complex plane which separates 
the poles
\begin{equation}\label{poles_b_jl}
b_{jl}=\frac{-b_j-l}{\beta_j}\quad (j=1,\dots,m;\, l=0,1,2,\dots)
\end{equation}
of the Gamma function \(\Gamma(b_j+\beta_j s)\) to the left of \(\mathcal{L}\) 
and the poles
\begin{equation}\label{poles_a_ik}
a_{ik}=\frac{1-a_i+k}{\alpha_i} \quad (i=1,\dots,n;\, k=0,1,2,\dots)
\end{equation}
to the right of \(\mathcal{L}\), and has one of the following forms
\begin{itemize}
\item[(i)] \(\mathcal{L}=\mathcal{L}_{-\infty}\) is a left loop situated in a horizontal strip starting from \(-\infty-\mathrm{i}\epsilon\) and terminating at the point \(-\infty+\mathrm{i}\epsilon\) with an \(\epsilon>0\).
\item[(ii)] \(\mathcal{L}=\mathcal{L}_\infty\) is a right loop situated in a horizontal strip starting at the point \(+\infty-\mathrm{i}\epsilon\) and terminating at the point \(+\infty+\mathrm{i}\epsilon\) with an \(\epsilon>0\).
\item[(iii)] \(\mathcal{L}=\mathcal{L}_{\mathrm{i}\gamma\infty}\) is a vertical line starting at the point \(\gamma-\mathrm{i}\infty\) and terminating at the point \(\gamma+\mathrm{i}\infty\), where \(\gamma\in\mathbb{R}\). 
\end{itemize}

Actually in the general case \(\mathcal{H}^{mn}_{pq}\) is not necessarily the Mellin transform of the $H$ function \(H^{mn}_{pq}\), but this is the case, when \(a^*>0\) with \(a^*\) defined in~\eqref{astar_def}~\cite[Theorem 2.2]{KilbSaig04}. In our considerations the condition \(a^*>0\) is valid, so we are always allowed to say that \(H^{mn}_{pq}\) and \(\mathcal{H}^{mn}_{pq}\) form a Mellin transform pair. The Mellin transform is a central tool in fractional calculus. In our considerations it is useful via its connection to the $H$ functions. For other useful applications we refer to~\cite{KiriLuch13}. The Mellin convolution relation
\begin{equation}\label{mellin_convolution_relation}
\mathcal{M}\left(f\stackrel{\mathcal{M}}{*} g\right)(s)=\mathcal{M}(f)(s)\mathcal{M}(g)(s),
\end{equation}
where 
\begin{equation}\label{mellin_convolution}
\left( f\stackrel{\mathcal{M}}{*} g\right)(x)=\int_0^\infty f\left(\frac{x}{y}\right) g(y)\frac{\mathrm{d}y}{y}
\end{equation}
is the {\em Mellin convolution}, is a very useful property. In particular, it is immediate from the definition of the $H$ function that the family of $H$ functions is invariant under the Mellin convolution. The Mellin convolution relation plays the same role as the convolution theorem for the Fourier transform. 

Since in our case the numbers \(a_i,\alpha_i,b_i,\beta_i\) are real, the poles \(a_{ik}\) in~\eqref{poles_a_ik} and \(b_{jl}\) in~\eqref{poles_b_jl} lie on the real axis, so by Cauchy's integral theorem one can change the infinite contour \(\mathcal{L}\) from one of the cases (i)--(iii) to another provided the integral in~\eqref{mellin_barnes_integral} converges in these cases. We collect here some of the results. The conditions for the convergence depend on the parameters \((a_i,\beta_i)_{1,p},(b_i,\beta_i)_{1,q}\) and can be derived from the Stirling formula for the Gamma function.  For further details we refer to~\cite{KilbSaig04} and references therein.

In the analysis we use the following properties from Chapter 2 of~\cite{KilbSaig04}. 


\begin{lemma}\label{Fox_properties}
Properties of Fox H-functions:
\begin{list2}
\vspace{.1in}
\item[(i)] For \(\omega,c\in\mathbb{C}\) and \(\sigma>0\) there holds 
\begin{align*}
\frac{\mathrm{d}^k}{\mathrm{d} z^k}\left\{z^\omega H^{mn}_{pq}\left( cz^\sigma \big| \begin{smallmatrix}(a_i,\alpha_i)_{1,p}\\
(b_j,\beta_j)_{1,q} \end{smallmatrix} \right)\right\} &= z^{\omega-k}H^{m, n+1}_{p+1, q+1}\left( cz^\sigma \big| \begin{smallmatrix} (-\omega,\sigma), &(a_i,\alpha_i)_{1,p}\\
(b_j,\beta_j)_{1,q}, &(k-\omega,\sigma) \end{smallmatrix} \right)\\
&=(-1)^k z^{\omega-k}H^{m+1, n}_{p+1, q+1}\left( cz^\sigma \big| \begin{smallmatrix} (a_i,\alpha_i)_{1,p}, &(-\omega,\sigma)\\
(k-\omega,\sigma), &(b_j,\beta_j)_{1,q} \end{smallmatrix} \right).
\end{align*} 
\vspace{.2in}

\item[(ii)] For \(m\ge 1\) and \(p>n\) there holds
\[ H_{pq}^{mn}\left( z \big| \begin{smallmatrix}
                            (a_i,\alpha_i)_{1,p-1},&(b_1,\beta_1)\\
(b_j,\beta_j)_{1,q} &
                           \end{smallmatrix}
\right)=H_{p-1,q-1}^{m-1,n}\left( z \big| \begin{smallmatrix}
                            (a_i,\alpha_i)_{1,p-1}\\
(b_j,\beta_j)_{2,q}
                           \end{smallmatrix}
\right). \]
\vspace{.2in}

\item[(iii)] For \(n\ge 1\) and \(q>m\) there holds
\[ H_{pq}^{mn}\left( z \big| \begin{smallmatrix}
                            (a_1,\alpha_1), &(a_i,\alpha_i)_{2,p}\\
(b_j,\beta_j)_{1,q-1} & (a_1,\alpha_1)
                           \end{smallmatrix}
\right)=H_{p-1,q-1}^{m,n-1}\left( z \big| \begin{smallmatrix}
                            (a_i,\alpha_i)_{2,p}\\
(b_j,\beta_j)_{1,q-1}
                           \end{smallmatrix}
\right). \] 
\vspace{.2in}

\item[(iv)] \(H_{pq}^{mn}\left( z^{-1} \big| \begin{smallmatrix} (a_i,\alpha_i)_{1,p} \\
                                 (b_j,\beta_j)_{1,q}
                                \end{smallmatrix}
\right)=H_{qp}^{nm}\left( z \big| \begin{smallmatrix} 
(1-b_j,\beta_j)_{1,q}\\
(1-a_i,\alpha_i)_{1,p}                                  
                                \end{smallmatrix}
\right). 
\)
\vspace{.2in}

\item[(v)] 
\(z^\sigma H_{pq}^{mn}\left( z \big| \begin{smallmatrix}
                            (a_i,\alpha_i)_{1,p}\\
(b_j,\beta_j)_{1,q} &
                           \end{smallmatrix} \right)=H_{pq}^{mn}\left( z \big| \begin{smallmatrix}
                            (a_i+\sigma\alpha_i,\alpha_i)_{1,p}\\
(b_j+\sigma\beta_j,\beta_j)_{1,q} &
                           \end{smallmatrix}\right) \) for \(\sigma\in\mathbb{C}\).

\vspace{.1in}

\item[(vi)] 
\(H_{pq}^{mn}\left( z \big| \begin{smallmatrix}
                            (a_i,\alpha_i)_{1,p}\\
(b_j,\beta_j)_{1,q} &
                           \end{smallmatrix} \right)=kH_{pq}^{mn}\left( z^k \big| \begin{smallmatrix}
                            (a_i,k\alpha_i)_{1,p}\\
(b_j,k\beta_j)_{1,q} &
                           \end{smallmatrix}\right) \) for \(k>0\).

\vspace{.1in}

\item[(vii)] For $b>0$ and $x>0$ there holds
\begin{align*}
&\int_0^\infty (xr)^{\omega}J_\eta(xr) H_{pq}^{mn}\left( br^\tau \big|\begin{smallmatrix}(a_i,\alpha_i)_{1,p}\\
(b_j,\beta_j)_{1,q} \end{smallmatrix} \right) \, \mathrm{d}r \\
&=\frac{2^\omega}{x}H_{p+2, q}^{m, n+1}\left( b2^\tau x^{-\tau} \big|\begin{smallmatrix} \l(1-\frac{\omega+1}{2}-\frac{\eta}{2}, \frac{\tau}{2}\r), &(a_i,\alpha_i)_{1,p}, &\l(1 -\frac{\omega+1}{2}+\frac{\eta}{2}, \frac{\tau}{2}\r)\\
(b_j,\beta_j)_{1,q} \end{smallmatrix} \right)
\end{align*}
provided the integral on the left hand side converges absolutely.

\vspace{.1in}
\end{list2}

\begin{proof}
The first six properties are straightforward calculations based on the Mellin-Barnes integral representation~\eqref{mellin_barnes_integral} of Fox $H$-functions. Indeed, property \((i)\) follows from~\eqref{mellin_barnes_integral} and the differentiation rule
\begin{align*}
\frac{\mathrm{d}^k}{\mathrm{d} z^k} z^{\omega-\sigma s}&=(\omega-\sigma s)(\omega-\sigma s-1)\cdots (\omega-\sigma s-k+1)z^{\omega-\sigma s-k}\\
&=z^{\omega-k}\frac{\Gamma(1+\omega-\sigma s)}{\Gamma(1+\omega-k-\sigma s)}z^{-\sigma s},
\end{align*}
where we used the property \(\Gamma(z+1)=z\Gamma(z)\) of the Gamma function.

Properties \((ii)\) and \((iii)\) follow from~\eqref{mellin_transform} by cancelling the common factors either from \(A(s)\) and \(C(s)\), or from \(B(s)\) and \(D(s)\).

Properties \((iv)\) and \((v)\) follow by the simple change of variables.

Property \((vii)\) follows by using the definition of the \(H\)-function~\eqref{mellin_barnes_integral}, changing the order of integration and using a known formula for the Mellin transform of the Bessel function \(J_\eta\). The conditions for the convergence follow from the asymptotics of \(J_\eta\) and the Stirling formula for the Gamma function. For details we refer to~\cite{KilbSaig04}, in particular Corollary 2.5.1. 
\end{proof}
\end{lemma}

As Theorem~\ref{thm_Bernstein} indicates, the Laplace transform plays a central role in our analysis. We will need the following result for the Laplace transform of the $H$ function~\cite[Section 2.5]{KilbSaig04}.

\begin{theorem}\label{thm_H_function_laplace}
Let \(a^*>0\) with \(a^*\) defined by~\eqref{astar_def}. Assume
\[
\min_{1\le j\le m}\left\lbrace\frac{\mathrm{Re}(b_j)}{\beta_j}\right\rbrace>-1.
\]
Then the Laplace transform of the $H$ function exists and there holds the relation
\begin{equation}\label{H_laplace_trans}
\left(\mathcal{L} H^{mn}_{pq}\left(x\Big| \begin{smallmatrix}
                           (a_i,\alpha_i)_{1,p}\\
(b_j,\beta_j)_{1,q}                 
                                           \end{smallmatrix}\right)\right)(p)=\frac{1}{p}H^{m,n+1}_{p+1,q}\left(\frac{1}{p}\Big|\begin{smallmatrix}
                           (0,1), &(a_i,\alpha_i)_{1,p}\\
(b_j,\beta_j)_{1,q} &                
                                           \end{smallmatrix}\right) 
\end{equation}
for \(p\in\mathbb{C}\) with \(\mathrm{Re}(p)>0\).
\end{theorem}

\subsection{Asymptotic behavior of the Fox 
$H$-functions}\label{section_Fox_asymptotics}
The asymptotic behavior of the $H$ functions play a central role in our analysis. The algebraic asymptotic expansions can be derived from the series expansions given by Cauchy's residue theorem and the results presented in~\cite{KilbSaig04}.  We introduce the following parameters
\begin{align}
a^*&=\sum_{i=1}^n\alpha_i-\sum_{i=n+1}^p \alpha_i+\sum_{j=1}^m\beta_j-\sum_{j=m+1}^q\beta_j,\label{astar_def}\\
D&=\sum_{j=1}^q\beta_j-\sum_{i=1}^p\alpha_i,\label{D_def}\\
\delta&=\prod_{i=1}^p\alpha_i^{-\alpha_i}\prod_{j=1}^q\beta_j^{\beta_j},\label{delta_def}\\
\mu&=\sum_{j=1}^m b_j-\sum_{i=1}^p a_i+\frac{p-q}{2}\label{mu_def}
\end{align}
for the decription of the results.

The Mellin transform \(\mathcal{H}^{mn}_{pq}\) given by~\eqref{mellin_transform} of the $H$-function given by~\eqref{mellin_barnes_integral} has the following asymptotic behavior~\cite[Chapter 1]{KilbSaig04}.

\begin{lemma}\label{mellin_transform_asymptotics}
For \(s=r+\mathrm{i}\rho\) there holds the estimates
\[
|\mathcal{H}^{mn}_{pq}(s)|\sim\left(\frac{\mathrm{e}}{|r|}\right)^{\mp D|r|}\delta^{\pm |r|}|r|^{\mathrm{Re}(\mu)},\quad r\to\pm\infty, 
\]
and
\[
|\mathcal{H}^{mn}_{pq}(s)|\sim |\rho|^{Dr+\mathrm{Re}(\mu)}\mathrm{e}^{\pi |\rho|a^*/2},\quad |\rho|\to\pm\infty, 
\]
uniformly in \(r\) on any bounded interval in \(\mathbb{R}\).
\end{lemma}

With this Lemma and the theory of residues one can prove the following result~\cite[Chapter 1]{KilbSaig04}.

\begin{theorem}\label{fox_series_expansion}
Let the parameters \(a^*,D\) and \(\delta\) be given by~\eqref{astar_def},~\eqref{D_def} and~\eqref{delta_def}. Define \(H^{mn}_{pq}(z)\) by~\eqref{mellin_barnes_integral} with the contour \(\mathcal{L}\) specified in the following cases.
\begin{itemize}
\item[(i)] Suppose that either \(D>0\) and \(z\neq 0\), or \(D=0\) and \(0<|z|<\delta\) hold. Then \(H^{mn}_{pq}(z)\) defines an
analytic function of \(z\) for \(\mathcal{L}=\mathcal{L}_{-\infty}\) and
\begin{equation}\label{Fox_series}
H^{mn}_{pq}(z)=\sum_{j=1}^m\sum_{l=0}^\infty 
\mathrm{Res}_{s=b_{jl}}\left(\mathcal{H}^{mn}_{pq}(s)z^{-s}\right),
\end{equation}
where \(b_{jl}\) are given in~\eqref{poles_b_jl}. 
\item[(ii)] Suppose that either \(D<0\) and \(z\neq 0\), or \(D=0\) and \(|z|>\delta\) hold. Then \(H^{mn}_{pq}(z)\) defines an analytic function of \(z\) for \(\mathcal{L}=\mathcal{L}_\infty\) and 
\begin{equation}\label{Fox_series2}
H^{mn}_{pq}(z)=-\sum_{i=1}^n\sum_{k=0}^\infty \mathrm{Res}_{s=a_{ik}}\left(\mathcal{H}^{mn}_{pq}(s)z^{-s}\right)
\end{equation}
\item[(iii)] If \(a^*>0\), then \(H^{mn}_{pq}(z)\) defines an analytic function of \(z\) for \(\mathcal{L}=\mathcal{L}_{\mathrm{i}\gamma\infty}\) in the sector \(|\arg z|<a^*\pi/2\).
\end{itemize}
\end{theorem}

\begin{remark}\label{rem_aD}
Note that the parameters \(a^*\) and \(D\) for the $H$ function $H^{21}_{23}$ appearing in the representation formula~\eqref{fund_solution} read \(a^*=2-\alpha\) and \(D=\beta-\alpha\), so by Theorem~\ref{fox_series_expansion} the $H$ function $H^{21}_{23}$ is well-defined for all \(\alpha\in(0,2)\) and \(\beta\in (0,2]\).
\end{remark}

The algebraic asymptotic behavior of \(H^{mn}_{pq}(z)\), as \(z\to 0\), follows 
immediately from~\eqref{Fox_series} in the case \(D\ge 0\) by 
calculating the residues. Similarly, the algebraic asymptotic behavior of \(H^{mn}_{pq}(z)\) as \(z\to\infty\) follows immediately from~\eqref{Fox_series2} in the case \(D\le 0\) by calculating the residues.

The asymptotic behavior of \(H^{mn}_{pq}(z)\) at infinity, when \(D\ge 0\) and \(\mathcal{L}=\mathcal{L}_\infty\) follow from Cauchy's theorem, which allows for \(a^*>0\) to change the contour \(\mathcal{L}=\mathcal{L}_{-\infty}\) to \(\mathcal{L}_{\mathrm{i}\gamma\infty}\) and to continue \(H^{mn}_{pq}(z)\) analytically to the sector \(|\arg z|<a^*\pi/2\), see~\cite[Section 1.5]{KilbSaig04}.


\begin{theorem}[The algebraic asymptotic behavior at infinity]\label{Braa36}
Let either \(D\le 0\), or \(D>0\) and \(a^*>0\) with \(a^*\) and \(D\) given by~\eqref{astar_def} and~\eqref{D_def}. The asymptotic expansion at infinity of the $H$-function \(H^{mn}_{pq}(z)\) defined by~\eqref{mellin_barnes_integral} is given by the series expansion~\eqref{Fox_series2}
in the sector \(|\arg z|<a^*\pi/2\).
\end{theorem}

When \(D<0\) and \(a^*>0\), the algebraic asymptotic expansion of \(H^{mn}_{pq}(z)\) near zero can be argued similarly as above appealing to analytic continuation and the Cauchy theorem. The asymptotics is given again by the series~\eqref{Fox_series} in the sector \(|\arg z|<a^*\pi/2\).

In some cases of parameters the $H$ function has exponential asymptotic behavior, which is more involved and we do not discuss it here. For details we refer to~\cite{Braa36} and 
~\cite[Sections 1.6 and 1.7]{KilbSaig04}. We just give here the result we need.

\begin{theorem}[The exponential asymptotic behavior at infinity]\label{thm_exp_asymptotics}
Let \(D>0\) and \(a^*>0\) with \(a^*\) and \(D\) given by~\eqref{astar_def} and~\eqref{D_def}. The $H$-function \(H^{m0}_{pm}(z)\) defined by~\eqref{mellin_barnes_integral} has for some positive constant \(c\) the following bound 
\begin{equation}\label{exp_bound_infinity}
|H^{m0}_{pm}(z)|\lesssim \exp(-cz^{1/D}),\quad z\to\infty, 
\end{equation}
on the sector \(|\arg z|<D\pi/2\) provided \(\mu\in (-1,0)\) with \(\mu\) given by~\eqref{mu_def}.
\end{theorem}

\subsection{The Mittag-Leffler function}\label{sec:ML_function}

An important special function in the fractional calculus is the 
Mittag-Leffler function
\begin{equation}\label{mittag_asymp1}
E_{\alpha}(z)=\sum_{k=0}^\infty\frac{z^k}{\Gamma(1+\alpha k)},\quad z\in\mathbb{C},
\end{equation}
which may viewed as the generalization of the exponential function, since clearly \(E_1(z)=\exp(z)\). Using Cauchy's residue theorem one can show that~\cite[Formulae (7.79) and (7.80)]{Mari83}
\begin{equation}\label{ML_fox}
E_\alpha(-z)=\frac{1}{2\pi\mathrm{i}}\int_{\mathcal{L}_{-\infty}}\frac{\Gamma(s)\Gamma(1-s)}{\Gamma(1-\alpha s)}z^{-s}\mathrm{d}s=H^{11}_{12}\left(z \big| \begin{smallmatrix} 
(0,1)& \\ (0,1), &(0,\alpha)
\end{smallmatrix}\right),
\end{equation}
Hence the Mittag-Leffler function is also a special case of the $H$-function. If we further specialize \(\alpha=1\) in~\eqref{ML_fox} and use the property $(ii)$ of Lemma~\ref{Fox_properties}, we obtain the exponential function as a special case of the $H$-function,
\begin{equation}\label{exp_fox}
\exp(-z)=H^{10}_{01}\left(z \big| \begin{smallmatrix} 
- \\ (0,1)
\end{smallmatrix}\right),
\end{equation}

The function $E_{\alpha}$ has the asymptotic behavior
\begin{equation}\label{ML2_estimate}
E_{\alpha}(-x) \sim \frac{1}{\Gamma(1-\alpha)}x^{-1}, \quad x\to\infty,\quad 1\neq\alpha\in (0,2).
\end{equation}

The asymptotic behavior~\eqref{ML2_estimate} follows from the integral 
representation
\[
E_{\alpha}(z)=\frac{1}{2\pi 
\mathrm{i}}\int_{\mathcal{C}}\frac{t^{\alpha-1}e^t}{t^\alpha-z}dt,
\]
where \(\mathcal{C}\) is an infinite contour in the complex plane. For details 
we refer to~\cite[Chapter 18]{Erdelyi53} and~\cite[Chapter 1]{Podl99}. Alternatively, one  can use the 
connection~\eqref{ML_fox} to the \(H\)-functions and use Theorem~\ref{Braa36}.

\def\cprime{$'$} \def\cprime{$'$}

\bibliographystyle{plain}

$\mbox{}$

\noindent {\footnotesize {\bf Jukka Kemppainen}, Applied and Computational Mathematics,
P.O. Box 8000,
90014 University of Oulu,
Finland
e-mail: Jukka.T.Kemppainen@oulu.fi





}

\end{document}